\documentstyle[fleqn,12pt]{article}
\begin{document}
\pagenumbering{arabic}
\setcounter{page}{1}
\pagestyle{plain}
\baselineskip=16pt

\thispagestyle{empty}
\rightline{MSUMB 98-02, May 1998} 
\vspace{1.4cm}

\begin{center}
{\Large\bf Differential Geometry of the $q$-superplane}
\end{center}

\vspace{1cm}
\begin{center} Salih \c Celik \footnote{E-mail: scelik@fened.msu.edu.tr}
\footnote{New E-mail address: sacelik@yildiz.edu.tr}
 
Mimar Sinan University, Department of Mathematics, \\
80690 Besiktas, Istanbul, TURKEY. \end{center}

\vspace{2cm}
{\bf Abstract}. 
Hopf algebra structures on the extended $q$-superplane and its differential 
algebra are defined. An algebra of forms which is obtained from the 
generators of the extended $q$-superplane is introduced and its Hopf algebra 
structure is given. 

\vfill\eject
\noindent
{\bf 1. Introduction}

\noindent
Differential geometry of Lie groups plays an important role in the mathematical 
modelling of physics theories. A class of noncommutative Hopf algebra has been 
found in the discussions of integrable systems. These Hopf algebras are 
$q$-deformed functions algebras of classical groups and this structure is 
called a quantum group [1]. The quantum group can also be regarded as a 
generalization of the notion of a group [2]. Thus it is also attractive to 
generalize the corresponding notions of differential geometry. Mathematical 
aspects of such a generalization are promising. More recently it has been 
suggested that the zero branes in matrix theory [3] should be identified with 
supercoordinates in noncommutative geometry [4]. 

Noncommutative geometry [4] has started to play an important role 
in different fields of mathematical physics for the last few years. 
The basic structure giving a direction to the noncommutative geometry 
is a differential calculus on an associative algebra. The noncommutative 
differential geometry of quantum groups was introduced by 
Woronowicz in ref. 5. In this approach the quantum group is taken as the 
basic noncommutative space and the differential calculus on the group 
is deduced from the properties of the group. The other approach, 
initiated by Wess and Zumino [6], succeed Manin's emphasis [7] on the 
quantum spaces as the primary objects, differential forms are defined 
in terms of noncommuting (quantum) coordinates, and the differential 
and algebraic properties of quantum groups acting on these spaces are 
obtained from the properties of the spaces. The natural extension of 
their scheme to superspace [8] was introduced by Soni in ref. 9. 

The quantum superplane is the simplest example of a noncommutative 
superspace. We have investigated the noncommutative geometry of the 
quantum superplane. In section 2 we introduce two noncommutative differential 
calculi on the $q$-superplane. One of them is quite different from 
the calculus described in ref. 9, where $GL_q(1\vert 1)$ covariance was 
assumed. The graded Hopf algebra structures of the extended $q$-superplane 
and these super-calculi are given in section 3. In the following section we 
introduce two forms from the differential algebra and also give the 
graded Hopf algebra structure of the obtained algebra of forms. 

\noindent
{\bf 2. Differential calculi on $q$-superplane} 

\noindent
Let us begin with the Manin superplane. The quantum superplane is defined 
as an associative algebra whose the even coordinate $x$ and the odd 
(grassmann) coordinate $\theta$ satisfy 
$$ x \theta - q \theta x = 0 \qquad \theta^2 = 0 \eqno(1) $$
where $q$ is a nonzero complex parameter. The algebra of 
$q$-polynomials will be called the algebra of functions on the quantum 
2-dimensional supervector space (superplane) and will be denoted by 
${\cal A}$. 

In order to establish a noncommutative differential calculus on the quantum 
superplane, we assume that the commutation relations between the 
coordinates and their differentials are in the following form 
$$x ~{\sf d}x = A {\sf d}x ~x $$
$$x ~{\sf d}\theta = F_{11} {\sf d}\theta ~x + F_{12} {\sf d}x ~\theta 
  \eqno(2)$$
$$\theta ~{\sf d}x = F_{21} {\sf d}x ~\theta + F_{22} {\sf d}\theta ~x $$
$$\theta ~{\sf d}\theta = B {\sf d}\theta ~\theta. $$

The coefficients $A, B$ and $F_{ij}$ will be determined in terms of the 
complex deformation parameter $q$. To find them we shall use the 
consistency of calculus. We first note that the properties of the 
exterior differential. 
The exterior differential {\sf d} is an operator which gives the mapping from 
the generators of ${\cal A}$ to the differentials: 
$${\sf d} : u \longrightarrow {\sf d}u \qquad u \in \{x,\theta\}. \eqno(3)$$
We demand that the exterior differential {\sf d} has to satisfy two 
properties: the nilpotency 
$${\sf d}^2 = 0 \eqno(4)$$
and the graded Leibniz rule 
$${\sf d}(f g) = ({\sf d} f) g + (- 1)^{\hat{f}} ({\sf d} g) \eqno(5)$$
where $\hat{f} = 0$ for even variables and $\hat{f} = 1$ for odd variables. 
From the consistency conditions 
$${\sf d}(x \theta - q \theta x) = 0 \qquad {\sf d}(\theta^2) = 0 $$
we find 
$$F_{11} + q F_{22} = q \qquad F_{12} + q F_{21} = - 1 \qquad B = 1.
\eqno(6\mbox{a})$$
Similarly, from 
$$(x \theta - q \theta x) {\sf d}x = 0 \qquad 
  (x \theta - q \theta x) {\sf d}\theta = 0 $$
one has 
$$F_{12} F_{22} = 0 \qquad (F_{11} - q A) F_{22} = 0. \eqno(6\mbox{b})$$
The system (6) has, at least, two solutions and we shall discuss them below. 

We now define the commutation relations between variables and their 
differentials in the following form 
$$Z^i {\sf d}Z^j = (- 1)^{\hat{i}(\hat{j} + 1)} C^{ji}{}_{kl} 
  {\sf d}Z^k Z^l \eqno(7)$$
where $C \in End({\cal C}\otimes {\cal C})$. Comparing (7) with (2) we obtain 
the general matrix C 
$$C  
 = \left(\matrix{ A &       0    &      0   & 0 \cr 
                  0 & - F_{21} & - F_{22} & 0 \cr
                  0 &    F_{12}  &   F_{11} & 0 \cr 
                  0 &       0    &      0   & 1 \cr }\right). \eqno(8)$$
In the language of matrix $C$, associative and consistency with the properties 
of {\sf d} requires that $C$ fulfill the following conditions: 
$$C_{12} C_{13} C_{23} = C_{23} C_{13} C_{12} \qquad 
  \hat{C}_{12} \hat{C}_{23} \hat{C}_{12} = 
   \hat{C}_{23} \hat{C}_{12} \hat{C}_{23} \eqno(9)$$
where $C_{12} = C \otimes I$, etc., $\hat{C} = P C$ and $P$ is 
the superpermutation matrix. The general matrix $\hat{C}$ may have, 
at least, one of two distinct forms: 
$$\hat{C}_I = \left(\matrix{   
    p & 0       & 0       & 0 \cr 
    0 & 0       & pq      & 0 \cr
    0 & q^{-1}  & p - 1   & 0 \cr 
    0 & 0       & 0       & 1  \cr }\right) \qquad  F_{12} = 0 
\eqno(10\mbox{a})$$
and 
$$\hat{C}_{II} = \left(\matrix{   
    s & 0       & 0 & 0 \cr 
    0 & r       & q & 0 \cr
    0 & qr - 1  & 0 & 0 \cr 
    0 & 0       & 0 & 1  \cr }\right) \qquad  F_{22} = 0 
\eqno(10\mbox{b})$$
where $p, r, s \in {\cal C}$ are free parameters. Similar matrices are found 
in [10] to obtain differential calculi on the quantum plane. 

The matrix $\hat{C}_I$ satisfies all required conditions. If we set $s = qr$ 
then the matrix $\hat{C}_{II}$ also obeys all required conditions [11]. Each 
of these matrices leads to a family of differential calculi on the 
$q$-superplane. 

So we have the following commutation relations: for $\hat{C}_I$ 
$$x~ {\sf d}x = p {\sf d}x~ x \qquad 
  x~ {\sf d}\theta = pq {\sf d}\theta~ x \eqno(11\mbox{a})$$
$$\theta~ {\sf d}x = - q^{-1} {\sf d}x~ \theta + (1 - p) {\sf d}\theta~ x 
  \qquad \theta~ {\sf d}\theta = {\sf d}\theta~ \theta. $$
For $\hat{C}_{II}$ 
$$x~ {\sf d}x = s {\sf d}x~ x \qquad 
  x~ {\sf d}\theta = q {\sf d}\theta~ x + (qr - 1) {\sf d}x~ \theta  
  \eqno(11\mbox{b})$$
$$\theta~ {\sf d}x = - r {\sf d}x~ \theta 
  \qquad \theta~ {\sf d}\theta = {\sf d}\theta~ \theta. $$

In the case of Family I, it is easy to check that the differential structure 
is invariant under action of quantum supergroup $GL_q(1\vert 1)$ (see, e.g. 
[12]) if we take $p = q^{-2}$. Similarly one can see, in the case of 
Family II, that the differential structure is invariant under action of 
$GL_{q,r}(1\vert 1)$ (see, e.g. [13]) if we set $s = qr$. 

Applying the exterior differential {\sf d} to the first and second (or third) 
relations of (11) we obtain 
$$({\sf d}x)^2 = 0 \qquad {\sf d}x {\sf d}\theta = pq {\sf d}\theta {\sf d}x
  \eqno(12\mbox{a})$$
for Family I and 
$$({\sf d}x)^2 = 0 \qquad {\sf d}x {\sf d}\theta = 
  r^{-1} {\sf d}\theta {\sf d}x  \eqno(12\mbox{b})$$
for Family II. 

A differential algebra on the an associative algebra ${\cal B}$ is a 
$z_2$-graded associative algebra $\Gamma$ equipped with an operator {\sf d} 
that has the properties (3)-(5). Furthermore, the algebra $\Gamma$ has to be 
generated by $\Gamma^0 \cup \Gamma^1 \cup \Gamma^2$, where $\Gamma^0$ is 
isomorphic to ${\cal B}$. For ${\cal B}$ we write ${\cal A}$. Let us denote 
the algebra (as a matter of fact the module) generated by ${\sf d}x$ and 
${\sf d}\theta$ with the relations (11) by $\Gamma^1$, where $\Gamma^1$ is 
isomorphic to ${\sf d}{\cal A}$, and the algebra (12) by $\Gamma^2$. Let 
$\Gamma$ be the quoitent algebra of the free associative algebra on the set 
$\{x,\theta,{\sf d}x,{\sf d}\theta\}$ modulo the ideal $J$ that is generated 
by the relations (1), (11) and (12). 

In the following section we shall show that the algebra ${\cal A}$ 
($q$-superplane), the algebra $\Gamma^1$ and also $\Gamma^2$ are all the 
graded Hopf algebras and so is the algebra $\Gamma$. 

\noindent
{\bf 3. Hopf algebra structures}

\noindent
A Hopf algebra structure on the quantum plane was introduced in [14]. 
In this section we introduce a graded Hopf algebra structure on the 
algebra ${\cal A}$ (i.e. on the $q$-superplane) and give the natural 
extension on $\Gamma$. 

\noindent
{\it 3.1. A Hopf algebra structure on ${\cal A}$} 

\noindent
We know, from section 1, that the quantum superplane, ${\cal A}$, is an 
associative algebra over a field $k$ generated by two elements 
$x$, $\theta$ obeying the relations (1). We can now define a coproduct and a 
counit on the algebra ${\cal A}$ as follows. 

The coproduct $\Delta: {\cal A} \longrightarrow {\cal A} \otimes {\cal A}$ 
is defined by 
$$\Delta(x) = x \otimes x $$
$$\Delta(\theta) = \theta \otimes x + x \otimes \theta \eqno(13)$$
$$\Delta(1) = 1 \otimes 1. $$
The counit $\epsilon: {\cal A} \longrightarrow {\cal C}$ is given by 
$$\epsilon(x) = 1 \qquad \epsilon(\theta) = 0. \eqno(14)$$

The algebra ${\cal A}$ with the coproduct and the counit has a structrure of 
bi-algebra. One extends the algebra ${\cal A}$ by including inverse of $x$ 
which obeys 
$$x x^{-1} = 1 = x^{-1} x.$$
If we extend the algebra ${\cal A}$ by adding the inverse of $x$ 
then the algebra ${\cal A}$ admits a coinverse (antipode) 
$S: {\cal A} \longrightarrow {\cal A}$ defined by 
$$S(x) = x^{-1} \qquad S(\theta) = - x^{-1} \theta x^{-1}. \eqno(15)$$
The coinvers has the properties of an inverse and we have $S^2 = 1$. Indeed, 
$$S^{-1}(x) = S(x) \qquad S^{-1}(\theta) = S(\theta).$$
Note that 
$$\Delta(x^{-1}) = x^{-1} \otimes x^{-1}.$$
It is not difficult to verify the following properties of 
co-structures: 
$$(\Delta \otimes \mbox{id}) \circ \Delta = (\mbox{id} \otimes \Delta) \circ \Delta 
\eqno(16)$$
$$\mu \circ (\epsilon \otimes \mbox{id}) \circ \Delta 
  = \mu' \circ (\mbox{id} \otimes \epsilon) \circ \Delta \eqno(17)$$
$$m \circ (S \otimes \mbox{id}) \circ \Delta = \epsilon 
  = m \circ (\mbox{id} \otimes S) \circ \Delta \eqno(18)$$
where id denotes the identity mapping, 
$$\mu : {\cal C} \otimes {\cal A} \longrightarrow {\cal A} \qquad 
  \mu' : {\cal A} \otimes {\cal C} \longrightarrow {\cal A} $$
are the canonical isomorphisms, defined by 
$$\mu(k \otimes u) = ku = \mu'(u \otimes k) \qquad \forall u \in {\cal A} 
  \quad \forall k \in {\cal C} $$
and $m$ is the multiplication map 
$$m : {\cal A} \otimes {\cal A} \longrightarrow {\cal A} \qquad 
  m(u \otimes v) = uv. \eqno(19)$$
The multiplication in ${\cal A} \otimes {\cal A}$ follows the rule 
$$(A \otimes B) (C \otimes D) = (-1)^{\hat{B} \hat{C}} AC \otimes BD. 
  \eqno(20)$$

The coproduct, counit and coinverse which are specified above supply 
the algebra ${\cal A}$ with a graded Hopf algebra structure. 

\noindent
{\it 3.2. A Hopf algebra structure on $\Gamma$} 

\noindent
We first note that consistency of a differential calculus with commutation 
relations (1) means that the algebra $\Gamma$ is a graded associative algebra 
generated by the elements of the set $\{x,\theta,{\sf d}x,{\sf d}\theta\}$. 

Since the algebra $\Gamma$ is generated by the generators set 
$\{x,\theta,{\sf d}x,{\sf d}\theta\}$ we must only describe the actions of 
co-maps on the subset $\{{\sf d}x,{\sf d}\theta\}$. To denote the coproduct, 
counit and coinverse which will be defined on the algebra $\Gamma$ with those 
of ${\cal A}$ may be inadvisable. For this reason, we shall denote them with 
a different notation. To this end we define a map 
$\hat{\Delta}_R : \Gamma \longrightarrow \Gamma \otimes {\cal A}$ such that 
$$\hat{\Delta}_R \circ {\sf d} = ({\sf d} \otimes \mbox{id}) \circ \Delta. \eqno(21)$$
Thus we have 
$$\hat{\Delta}_R({\sf d}x) = {\sf d}x \otimes x $$
$$\hat{\Delta}_R({\sf d}\theta) = {\sf d}\theta \otimes x + 
 {\sf d}x \otimes \theta. \eqno(22)$$
We now define a map $\phi_R$ as follows: 
$$\phi_R(u_1 {\sf d}v_1 + {\sf d}v_2 u_2) = 
  \Delta(u_1) \hat{\Delta}_R({\sf d}v_1) + \hat{\Delta}_R({\sf d}v_2) 
  \Delta(u_2). \eqno(23)$$
Then it can be checked that the map $\phi_R$ leaves invariant the relations 
(11) and (12). One can also check that the following identities are satisfied: 
$$(\phi_R \otimes \mbox{id}) \circ \phi_R = 
  (\mbox{id} \otimes \Delta) \circ \phi_R 
  \qquad (\mbox{id} \otimes \epsilon) \circ \phi_R = \mbox{id}. \eqno(24)$$
But we do not have a coproduct for the differential algebra because the map 
$\hat{\Delta}_R$ does not gives an analog for the derivation property (5), 
yet. So we define another map 
$\hat{\Delta}_L : \Gamma \longrightarrow {\cal A} \otimes \Gamma$ such that 
$$\hat{\Delta}_L \circ {\sf d} = (\tau \otimes {\sf d}) \circ \Delta \eqno(25)$$
and a map $\phi_L$ with again (23) by replacing $L$ with $R$. 
Here $\tau$ is a map which gives $\tau(a) = (-1)^{\hat{a}} a$. 
The map $\phi_L$ also leaves invariant the relations (11) and (12), 
and the following identities are satisfied: 
$$(\mbox{id} \otimes \phi_L) \circ \phi_L = 
  (\Delta \otimes \mbox{id}) \circ \phi_L 
  \qquad (\epsilon \otimes \mbox{id}) \circ \phi_L = \mbox{id}. \eqno(26)$$

Let us define the map $\hat{\Delta}$ as 
$$\hat{\Delta} = \phi_R + \phi_L \eqno(27)$$
which will allow us to define the coproduct of the differential algebra. 
We denote the restriction of $\hat{\Delta}$ to the algebra ${\cal A}$ by 
$\Delta$ and the extension of $\Delta$ to the differential algebra $\Gamma$ by 
$\hat{\Delta}$: 
$$\hat{\Delta}\vert_{\cal A} = \Delta \qquad 
  \Delta\vert_\Gamma = \hat{\Delta}. \eqno(28)$$
It is possible to interpret the first relation in (28) as the definition of 
$\hat{\Delta}$ and (27) as the definition of $\hat{\Delta}$ on differentials. 

One can see that $\hat{\Delta}$ is a linear map and a homomorphism. 
In fact, for example, 
$$\hat{\Delta}(x~ {\sf d}x) = (\phi_R + \phi_L)(x~ {\sf d}x) = 
  \Delta(x) (\hat{\Delta}_R + \hat{\Delta}_L)({\sf d}x)$$
and with (28) 
$$\Delta(x) \hat{\Delta}({\sf d}x) = 
  \Delta(x) [\Delta(1)(\hat{\Delta}_R + \hat{\Delta}_L)({\sf d}x)].$$
Using the coassociativity of $\Delta$, eq. (16), we can also show the 
coassociativity of $\hat{\Delta}$ . So the map $\hat{\Delta}$ is a coproduct 
for the differential algebra $\Gamma$. 

Similarly, if we define a counit $\hat{\epsilon}$ for the differential 
algebra as 
$$\hat{\epsilon} \circ {\sf d} = {\sf d} \circ \epsilon = 0 \eqno(29)$$
and 
$$\hat{\epsilon}\vert_{\cal A} = \epsilon \qquad 
  \epsilon\vert_\Gamma = \hat{\epsilon}. \eqno(30)$$
one has 
$$\hat{\epsilon}({\sf d}x) = 0 \qquad \hat{\epsilon}({\sf d}\theta) = 0 
  \eqno(31)$$
where 
$$\hat{\epsilon}(u_1 {\sf d}v_1 + {\sf d}v_2 u_2) = 
  \epsilon(u_1) \hat{\epsilon}({\sf d}v_1) + 
  \hat{\epsilon}({\sf d}v_2) \epsilon(u_2). \eqno(32)$$
Here we used the fact that ${\sf d}(1) = 0$. 

The next step is to obtain a coinverse $\hat{S}$. For this, it suffices to 
define $\hat{S}$ such that 
$$\hat{S} \circ {\sf d} = {\sf d} \circ S \eqno(33)$$
and 
$$\hat{S}\vert_{\cal A} = S \qquad 
  S\vert_\Gamma = \hat{S} \eqno(34)$$
where 
$$\hat{S}(u_1 {\sf d}v_1 + {\sf d}v_2 u_2) = 
  \hat{S}({\sf d}v_1) S(u_1) + S(u_2) \hat{S}({\sf d}v_2). \eqno(35)$$
So the action of $\hat{S}$ on the generators ${\sf d}x$ and ${\sf d}\theta$ is 
as follows: 
$$\hat{S}({\sf d}x) = - x^{-1} ~{\sf d}x ~x^{-1} $$
$$\hat{S}({\sf d}\theta) = - x^{-1} ~{\sf d}\theta ~x^{-1} + 
  2 x^{-1} ~{\sf d}x ~x^{-1} \theta x^{-1}. \eqno(36)$$
Notice that it is easy to check that $\hat{\epsilon}$ and $\hat{S}$ leave 
invariant the relations (11) and (12). 

Consequently, we can say that the structure 
$(\Gamma, \hat{\Delta}, \hat{\epsilon}, \hat{S})$ is a graded Hopf algebra. 

\noindent
{\bf 4. Hopf algebra structure of forms on ${\cal A}$}

\noindent
In this section we shall define two forms using the generators of 
$\cal A$ and show that the algebra of forms is a graded Hopf algebra. 

If we call them $w$ and $u$ then one can define them as follows: 
$$w = {\sf d}x ~x^{-1} \qquad 
  u = {\sf d}\theta ~x^{-1} - {\sf d}x ~x^{-1} \theta x^{-1}. \eqno(37)$$
We denote the algebra of forms generated by two elements $w$ and $u$ by 
$\Omega$. The generators of the algebra $\Omega$ with the generators of 
$\cal A$ satisfy the following rules 

{\bf I}. 
$$x w = p w x \qquad \theta w = - w \theta + (1 - p) u x $$
$$x u = pq u x \qquad \theta u = pq u \theta. \eqno(38\mbox{a})$$

{\bf II}. 
$$x w = s w x \qquad \theta w = - qr w \theta $$
$$x u = q u x + q (qr - s) w \theta \qquad \theta u = q u \theta. 
 \eqno(38\mbox{b})$$
The commutation rules of the generators of $\Omega$ are 

{\bf I}. 
$$w^2 = 0 \qquad wu = uw. \eqno(39\mbox{a})$$

{\bf II}. 
$$w^2 = 0 \qquad wu = qrs^{-1} uw. \eqno(39\mbox{b})$$

We make the algebra $\Omega$ into a graded Hopf algebra with the following 
co-structures: the coproduct 
$\Delta: \Omega \longrightarrow \Omega \otimes \Omega$ is defined by 
$$\Delta(w) = w \otimes 1 + 1 \otimes w $$
$$\Delta(u) = u \otimes 1 + 1 \otimes u. \eqno(40) $$
The counit $\epsilon: \Omega \longrightarrow {\cal C}$ is given by 
$$\epsilon(w) = 0 \qquad \epsilon(u) = 0 \eqno(41)$$
and the coinverse $S: \Omega \longrightarrow \Omega$ is defined by 
$$S(w) = - w \qquad S(u) = - u. \eqno(42)$$
One can easy to check that (16)-(18) are satisfied. Note that the commutation 
relations (38) and (39) are compatible with $\Delta$, $\epsilon$ and $S$, 
in the sense that $\Delta(x w) = p \Delta(w x)$, $\Delta(w^2) = 0$ and so on. 

\noindent
{\bf Acknowledgement}

\noindent
This work was supported in part by TBTAK the Turkish Scientific and Technical 
Research Council. I would like to express my deep gratitude to the referees 
for critical comments and suggestions on the manuscript.

\end{document}